\RequirePackage{fix-cm}

\documentclass{elsarticle}
\usepackage{graphicx}
\usepackage{setspace}

\usepackage{float}
\usepackage{breakcites}
\usepackage{epsfig}
\usepackage{amssymb}
\usepackage{amsmath}
\usepackage{hyperref}
\usepackage{verbatim}
\usepackage{subcaption}
 \usepackage{natbib}
\usepackage{tikz}
\usepackage{color}

\usetikzlibrary{arrows}
\usepackage[bottom=0.95in,top=0.95 in]{geometry}
\newcommand\given[1][]{\:#1\vert\:}

\tikzset{cross/.style={cross out, draw=black, minimum size=2*(#1-\pgflinewidth), inner sep=0pt, outer sep=0pt},
cross/.default={5pt}}
\usetikzlibrary{shapes.misc}
\newcommand{\be}{\begin{eqnarray}}
\newcommand{\ee}{\end{eqnarray}}
\newcommand{\bea}{\begin{eqnarray*}}
\newcommand{\eea}{\end{eqnarray*}}

\usepackage{nicefrac}
\numberwithin{equation}{section}

\usepackage{pst-node,pst-plot}
\pstVerb{realtime srand}
\psset{plotpoints=50}

\begin{document}

\title{Approximating  Shepp's constants for the Slepian process
}

\author[1]{Jack Noonan}
 \ead{NoonanJ1@cardiff.ac.uk}
\author[1]{Anatoly Zhigljavsky\corref{cor2}}% \fnref{fn1}}
 \ead{ZhigljavskyAA@cardiff.ac.uk}

 \address[1]{School of Mathematics, Cardiff University, Cardiff, CF24 4AG, UK}

\cortext[cor2]{Corresponding author}


\begin{abstract}
Slepian process $S(t)$ is  a stationary Gaussian process with zero mean and covariance
$
\mathbb{E} S(t)S(t')=\max\{0,1-|t-t'|\}\, .
$
For any $T\geq 0$ and real $h$, define
$F_T(h ) = {\rm Pr}\left\{\max_{t \in [0,T]} S(t) < h   \right\}
$
 and the constants
 $\Lambda(h) = -\lim_{T \to \infty} \frac1T   \log F_T(h)$ and $\lambda(h)=\exp\{-\Lambda(h) \}$; we will call them `Shepp's constants'.
The aim of the paper is construction of accurate  approximations for $F_T(h)$ and hence for the Shepp's constants.
We
demonstrate that at least some of the approximations are extremely accurate.

\end{abstract}

\maketitle
\section{Introduction}
\label{sec:prob-state}

Let $S(t)$,  $ t\in [0,T]$, be  a Gaussian process with mean 0 and covariance
\be\label{cov}
\mathbb{E} S(t)S(t')=\max\{0,1-|t-t'|\}\, .
\ee
This process is often called Slepian process.
For any real $h$ and $x <h$, define
\be
F_T(h \given x ) := {\rm Pr}\left\{\max_{t \in [0,T]} S(t) < h \; \big | \; S(0) = x  \right\}\, ;
 \label{eq:prob-S}
\ee
if $x \geq h$ we set $F_T(h \given x )=0$.
Assuming that $x $ has Gaussian distribution $N(0,1)$, and hence the stationarity of the process $S(t)$, we average $F_T(h \given x )$ and thus define
\be
\label{eq:F_T}
F_T(h ):= \int_{-\infty}^h F_T(h \given x ) \varphi(x ) dx   \, ,
\ee
where $\varphi(x)=(2\pi)^{-1/2} \exp\{-x^2/2\}$.

Key results on the boundary crossing probabilities for the Slepian process have been established by L.Shepp in~\cite{Shepp71}. In particular,
Shepp has derived an explicit formula for $F_T(h)$ with $T$ integer, see \eqref{shepp_form} below.  As this explicit formula is quite complicated, in (3.7) in the same paper, Shepp has conjectured the existence of the following constant (depending on $h$)
\be
\label{eq:Shepp_constant}
\Lambda(h) = -\lim_{T \to \infty} \frac1T   \log F_T(h)
\ee
and raised the question of constructing accurate approximations and bounds for this constant.

The importance of this constant is related to the asymptotic relation
\be
\label{eq:Shepp_constant1}
F_T(h) \simeq {\rm const} [\lambda(h)]^{T} \, \;\;{\rm \; as}\;\; T \to \infty\, ,
\ee
where $\lambda(h)=\exp\{-\Lambda(h) \}$.
We will call  $\Lambda(h)$ and $\lambda(h)$ `Shepp's constants'.

In this paper,
we are interested in deriving approximations for $F_T(h)$ in the form \eqref{eq:Shepp_constant1} and hence for the Shepp's constants. In formulation of approximations, we offer approximations for $F_T(h)$ for all $T>2$ and hence approximations for $\Lambda(h)$ and $\lambda(h)$. Note that computation of $F_T(h)$ for $T\leq 2$ is a relatively easy problem, see \cite{slepian1961first} for $T\leq 1$ and \cite{Shepp71} for $1<T\leq 2$.

In Section~\ref{sec:2} we derive several approximations for $F_T(h)$ and  $\lambda(h)$ and  provide
numerical results showing  that at least some of the derived approximation are extremely accurate.
In Section~\ref{sec:4} we compare the upper tail asymptotics for the Slepian process and  some other stationary Gaussian processes.
Section~\ref{sec:app} contains some minor technical details and
Section~\ref{sec:conc} delivers conclusions.

\section{Construction of approximations}
\label{sec:2}

\subsection{Existence of Shepp's constants and the approximations derived from general principles}\label{sec:21}

The fact  that the limit in \eqref{eq:Shepp_constant} exists and hence that $\Lambda(h)$  is properly defined for any $h$ has been proven in \cite{li2004lower}. The proof of existence of $\Lambda(h)$ is based on the inequalities
\be\label{simple_bounds}
 -\frac1{n+1} \log [F_{n}(h)] \leq \Lambda(h) \leq -\frac1n \log [F_{n}(h)]\;\;{\rm for \; any\; }  n=1,2, \ldots
\ee
The inequality in the rhs of \eqref{simple_bounds} follows directly from the infamous `Slepian inequality' established in \cite{slepian1962}; this inequality holds for any  Gaussian stationary process with non-negative correlation function. The inequality in the lhs of \eqref{simple_bounds} can be obtained by a simple extension of the arguments in  \cite[p.470]{slepian1962}; it holds for any  Gaussian stationary process which correlation function vanishes outside the interval $[-1,1]$. The inequalities \eqref{simple_bounds} are not sharp: in particular, for  $n=2$ and $h=0$, \eqref{simple_bounds} gives  $1.336<\Lambda(0)< 2.004$; see \cite[Remark 3]{molchan2012survival}. As follows from Tables 1 and 3, an accurate approximation for $\Lambda(0)$   is $\Lambda(0) \simeq 1.5972,$ where we claim all four decimal places are accurate.

\begin{figure}[!h]
\label{bounds}
\begin{center}
\begin{subfigure}{0.48\linewidth} \centering
     \includegraphics[scale=0.2]{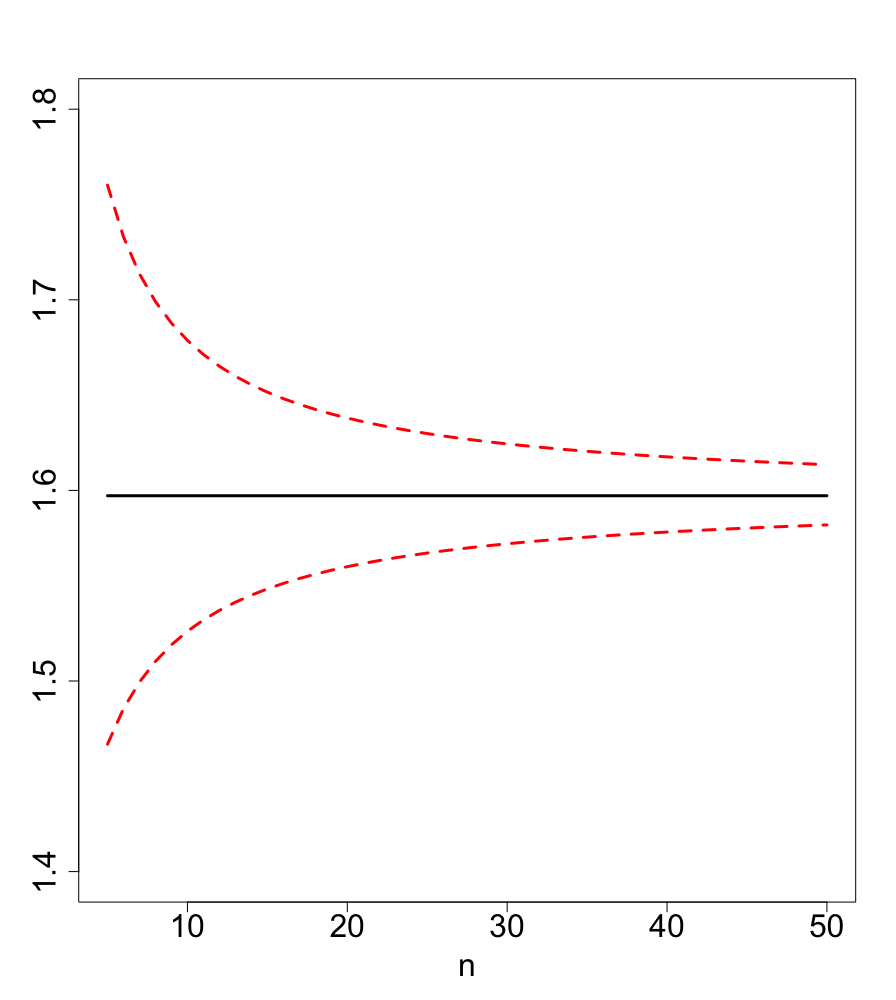}
     \caption{ $h=0$  } \label{fig:b1}
\end{subfigure}
\hspace{0.3cm}
 \begin{subfigure}{0.48\linewidth} \centering
     \includegraphics[scale=0.2]{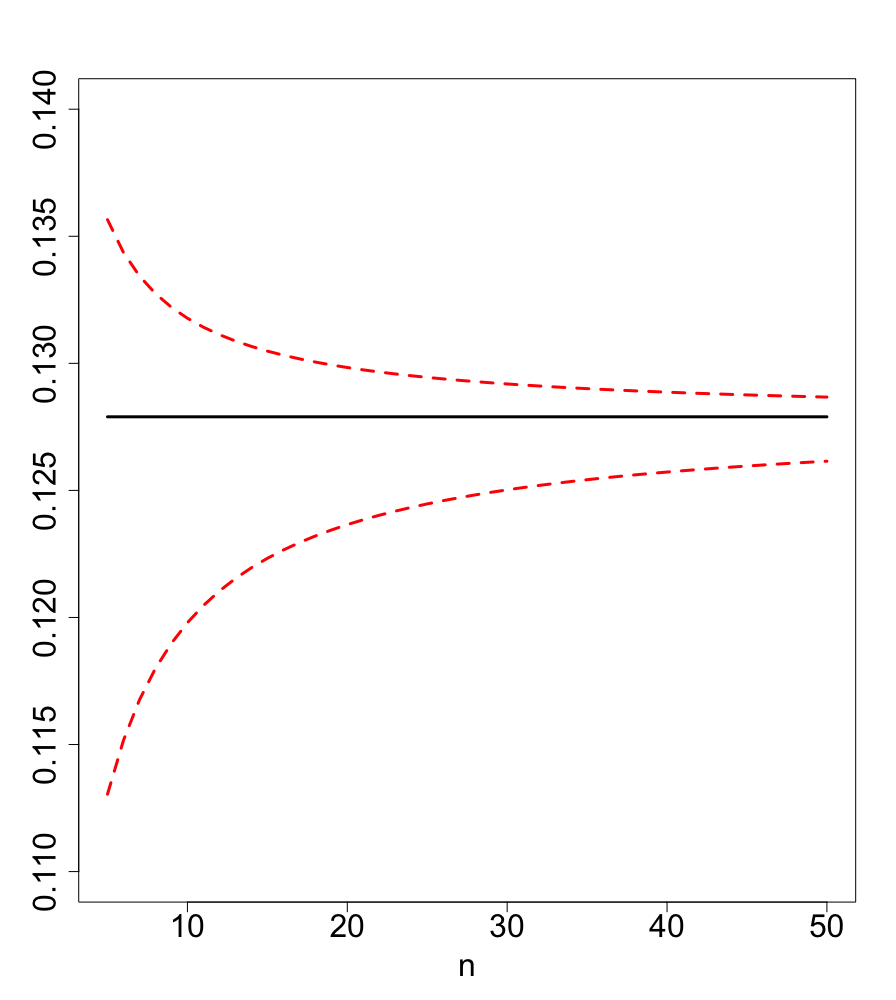}
     \caption{$h=2$} \label{fig:bounds}
\end{subfigure}
    \end{center}
    \vspace{-0.4cm}
    \caption{Lower and upper bounds \eqref{simple_bounds} (red dotted lines) for $\Lambda(h)$ (solid black line).  }
\end{figure}

If $n$ is not too small, the bounds \eqref{simple_bounds} are very difficult to compute. For small $h$, these bounds are not sharp even if $n$ is large, see Fig.~\ref{fig:b1}. The bounds improve as $h$ grows, see Fig.~\ref{fig:bounds}. It is not very clear how to use these bounds for construction of accurate approximations for $\Lambda(h)$. In particular, from
Fig.~\ref{fig:bounds} we  observe that  the upper bound of \eqref{simple_bounds} can be much closer to the true $\Lambda(h)$ than the lower bound.

One may apply general results shown in \cite{landau1970supremum,marcus1972sample}, see also  formula $(2.1.3)$ in \cite{adler2007random}, to approximate $F_T(h)$ for large $h$ but these results only show that $\lambda(h) \to 1$ as $h \to \infty$ and therefore are of no use here. A more useful tool, which can be used for approximating $\lambda(h)$, is connected to the following result of J.Pickands proved in
\cite{pickands1969upcrossing}. Assume that $\{\xi(t)\}$ is a stationary Gaussian random process with $\mathbb{E}\xi(t)=0$, $\mathbb{E}\xi^2(t)=1$ and covariance function \be
\label{eq:cov_f}
\rho(t)=\mathbb{E}\xi(0)\xi(t)=0=1-C|t|^\alpha + o(|t|^\alpha ) \; \mbox{ as $t \to 0$ }
\ee
and $\sup_{\epsilon\leq t \leq T}\rho(t) <1 ,$ $\forall \epsilon >0$. Then
\be
\label{eq:Pick}
{\rm Pr} \left\{ \sup_{0\leq t \leq T} \xi(t) \leq h\right\} = 1- T C^{1/\alpha} H_\alpha h^{2/\alpha-1} \varphi(h)  \left( 1+ o(1) \right) \;\;\;{\rm as}\; h \to \infty\, ,
\ee
where $H_\alpha$ is the so-called `Pickands constant'. By replacing $1-x$ with $e^{-x}$ ($x\to 0$) and  removing the term $\left( 1+ o(1) \right)$ in \eqref{eq:Pick} we obtain a general  approximation
\be \label{eq:Pick1}
{\rm Pr} \left\{ \sup_{0\leq t \leq T} \xi(t) \leq h\right\} \simeq \exp\{- T C^{1/\alpha} H_\alpha h^{2/\alpha-1} \varphi(h) \}\, .
\ee
As shown in \cite{harper2017pickands}, the value of the Pickands constant $H_\alpha$ is only known for $\alpha=1,2$ and hence  the approximation \eqref{eq:Pick1} can only be applied in these cases. When $\xi(t)$ is the Slepian process $S(t)$ with covariance function \eqref{cov} we have $\alpha=1$, $ H_1=1$ and $C=1$. Hence we obtain from  \eqref{eq:Pick1}
\bea
\mbox{\noindent{\bf Approximation 0:}}\;\;F_T(h) \simeq \exp(-h\varphi(h) T) \, ,\;\;\; \Lambda^{(0)}(h) = h\varphi(h)\, ,\;\; \lambda^{(0)}(h)= e^{-h\varphi(h) }\, .\;\;\;\;\;\;\;\;\;\;\;\;
\eea

Note that Approximation 0 can also be obtained as a Poisson clumping heuristic, see formula (D10g) in \cite{Aldous}.
If $h$ is not large, then Approximation 0 is quite poor, see Tables 1 and 2 and Figure~2.
%\ref{fig:lambda12}
% \ref{table1} and \ref{table2}
 For small and moderate values of $h$, the approximations derived below in this section are much superior to Approximation 0.

\subsection{Shepp's  formula for $F_n$}

%\subsubsection{Shepp's  formula }

The following formula is the result (2.15)  in  \cite{Shepp71}:
\begin{equation}\label{shepp_form}
F_n(h \big | x ) = \frac{1}{\varphi(x)} \int_{D_x} \det|\varphi(y_i - y_{j+1} + h){|^n} \hspace{-0.3cm}|_{i,j=0} \, dy_2\ldots dy_{n+1},
\end{equation}
where $T=n$ is a positive integer,
$
D_x =  \{y_2, \dots , y_{T+1} \given h-x < y_2 < y_3 < \ldots < y_{n+1}   \}
$, $y_0= 0, y_1=h-x.$ L.Shepp in \cite{Shepp71} has also derived explicit formulas for $F_T(h \big | x )$ with non-integral $T>0$ but these formulas are more complicated and are realistically  applicable only for small  $T$ (say, $T\leq 3$).

From \eqref{shepp_form} we straightforwardly obtain
\be\label{shepp_form1}
F_1(h \big | x ) &=& \Phi(h)-\frac{\varphi(h)}{\varphi(x)} \Phi(x)\, ,\\
F_1(h  ) &=& \int_{-\infty}^h F_1(h \big | x )\varphi(x) dx = \Phi^2(h)-\varphi(h) \big[h\Phi(h)
+ \varphi(h) \big]\, ,\label{eq:F1}
\ee
where $\Phi(x) = \int_{-\infty}^x \varphi(t) dt$.
Derivation of explicit formulas for $F_T(h \big | x )$ and $ F_T(h  )$ with $T\leq 1$ is relatively easy as the process $S(t)$ is conditionally Markovian in the interval $[0,1]$, see \cite{Mehr}.
Formula~\eqref{shepp_form1} has been first derived in \cite{slepian1961first}.

In what follows, $F_2(h  )$  also plays a very  important role.
Using \eqref{shepp_form} and changing the order of integration where suitable, $F_2(h)$ can be expressed through a one-dimensional integral as follows:
%\begin{eqnarray}
%F_2(h) &=& \Phi(h)^3 + \varphi(h)^2\Phi(h) + \frac{\varphi(h)^2}{2}\left[ (h^2-1)\Phi(h)+h\varphi(h) \right] + \int_{-\infty}^{h}\Phi(y)^2\varphi(2h-y)dy \nonumber\\
%&-& 2\varphi(h)\Phi(h)[h\Phi(h)+\varphi(h)] - \frac{1}{\sqrt{2}}\int_{0}^{\infty}\Phi(h-y)\varphi(\sqrt{2}h)\left[ \Phi(\sqrt{2}y)-1/2  \right]dy. \label{eq:F2}
%\end{eqnarray}

\begin{eqnarray}
F_2(h) &=& \Phi^3(h) + \varphi^2(h)\Phi(h) + \frac{\varphi^2(h)}{2}\left[ (h^2-1)\Phi(h)+h\varphi(h) \right] + \int_0^{\infty}\Phi^2(h-y)\varphi(h+y)dy \nonumber\\
&-& 2\varphi(h)\Phi(h)[h\Phi(h)+\varphi(h)] - \frac{1}{\sqrt{2}}\int_{0}^{\infty}\Phi(h-y)\varphi(\sqrt{2}h)\left[ \Phi(\sqrt{2}y)-1/2  \right]dy. \label{eq:F2}
\end{eqnarray}
This expression can be  approximated as shown in Appendix; see \eqref{Corrected_Diffusion_explicit2}.

\subsection{An alternative  representation of the Shepp's  formula \eqref{shepp_form}}

Let $T=n$ be a positive integer, $y_0= 0, y_1=h-x.$ For $i=0,1, \ldots, n$ we set $s_i=h+y_i-y_{i+1}$ with $s_0=x$.
It follows from Shepp's proof of \eqref{shepp_form} that  $s_0, s_1, \ldots, s_n$ have the meaning of the values of the process $S(t)$ at the times $t=0,1, \ldots, n$:
$S(i)=s_i$ ($i=0,1, \ldots, n$). The range of the variables $s_i$ is $(-\infty,h)$. The variables $y_1, \ldots, y_{n+1}$ are expressed via $s_0, \ldots, s_{n}$ by $y_k=kh-s_0-s_1-\ldots-s_{k-1}$
($k=1, \ldots, n+1$) with $y_0=0$. Changing the variables in \eqref{shepp_form}, we obtain 
\begin{equation}\label{shepp_form2}
F_n(h \big | x ) = \frac{1}{\varphi(x)} \int_{-\infty}^h \ldots \int_{-\infty}^h  \det|\varphi(s_i + a_{i,j}){|^n} \hspace{-0.3cm}|_{i,j=0} \, ds_1\ldots ds_{n}\,,
\end{equation}
where
\bea
a_{i,j}=y_{i+1} \!- \!y_{j+1}=\left\{ \begin{array}{cl}
 0 & \text{for } i=j\, \\
     (i-j)h \!- \!s_{j+1}\!-\!\ldots\!- \!s_{i+1} & \text{for } i>j\,  \\
       (i-j)h+  s_{i+1}+\ldots+ s_{j}  & \text{for } i<j\, .

       \end{array} \right.
\eea
Expression \eqref{shepp_form2} for the probability $F_n(h \big | x )$ implies that the function
\begin{equation}\label{shepp_form7}
p(s_0,s_1,\ldots s_{n})= \frac{1}{\varphi(s_0) F_n(h \big | s_0 )}  \det|\varphi(s_i +a_{i,j}){|^n} \hspace{-0.3cm}|_{i,j=0} \, .
\end{equation}
is the joint probability density function for the values $S(0),S(1), \ldots, S(n)$ under the condition  $S(t)<h$ for all $t\in [0,n]$.

Since $s_n$ is the value of $S(n)$, the formula \eqref{shepp_form7} also shows  the transition density from $s_0=x$ to $s_n$ conditionally
$S(t)<h$ for all $t\in [0,n]$:
 \begin{equation}\label{shepp_form9}
p^{(n)}_h(x\to s_n ) = \frac{1}{\varphi(x)} \int_{-\infty}^h \ldots \int_{-\infty}^h  \det|\varphi(s_i + a_{i,j}){|^n} \hspace{-0.3cm}|_{i,j=0} \, ds_1\ldots ds_{n-1}\, .
\end{equation}
For this transition density, $ \int_{-\infty}^h   p^{(n)}_h(x\to z ) dz = F_n(h \big | x )$.

\subsection{Approximating $\lambda(h)$ through eigenvalues of integral operators}\label{GL_approx}

\subsubsection{One-step transition}
\label{sec:one}

In the case $n=1$ we obtain from \eqref{shepp_form9}:
 \begin{equation}\label{shepp_form3}
p^{(1)}_h(x\to z ) = \frac{1}{\varphi(x)}   \det  \left(
                                                      \begin{array}{cc}
                                                       \varphi(x)  & \varphi(x\!-\!h\!+\!z) \\
                                                        \varphi(h) & \varphi(z) \\
                                                      \end{array}
                                                    \right)= \varphi(z)\left[1 - e^{-(h-z)(h-x)} \right]
\end{equation}
 with $z=s_1<h$.

Let $\lambda_1(h)$ be the largest eigenvalue of the
 the integral operator with kernel (\ref{shepp_form3}):
 \bea
\lambda_1(h) p(z) = \int_{-\infty}^{h}p(x)p^{(1)}_h(x\to z ) dx, \text{  } z<h\, , \label{2.27}
\eea
where  eigenfunction $p(x)$ is some probability density on $(-\infty,h]$.
The Ruelle-Krasnoselskii-Perron-Frobenius theory of bounded linear positive operators (see e.g. Theorem XIII.43 in \cite{ReedSimon}) implies
 that the maximum eigenvalue $\lambda$ of the operator with kernel $K(x,z)=p^{(1)}_h(x\to z )$ is simple, real and positive
 and the   eigenfunction $p(x)$ can be chosen as  a probability density.

 Similarly to what we have done  below in  Section~\ref{sec:two_step}, we can suggest computing good numerical approximations  to ${\lambda}_1(h)$ using Gauss-Legendre quadrature formulas.   However, we suggest to use (4.15) from~\cite{NoonZhig} instead; this  helps us to  obtain  the following  simple but rather accurate approximation to
 $\lambda_1(h)$:
$$
 \hat{\lambda}_1(h)  =\Phi(h) + {\varphi(h)}/{h} -\! { \varphi (h) [ \varphi (h) +h \Phi(h)] t
  }/
 { \left[\Phi(h)- e^{-h^2/2}/2 \right] } \, . \;\;
$$

\bea
\mbox{\noindent{\bf Approximation 1:}}\;\;F_T(h) \simeq  F_1(h) \left[ {\hat{\lambda}_1(h)} \, \right]^{T-1} \, \;\mbox{($T\geq 1$);} \;\;\Lambda^{(1)}(h) = - \log \hat{\lambda}_1(h) \, , \; \lambda^{(1)}(h)= {\hat{\lambda}_1(h)}\, .\;\;\;
\eea

%
%\noindent{\bf Approximation 1.}
%\bea
%F_T(h) \simeq  F_1(h) \left[ \hat{\lambda}_1(h)\right]^{T-1} \, \;\mbox{($T\geq 1$);}\;\;\; \Lambda^{(1)}(h) = -\log \hat{\lambda}_1(h) \, , \;\; \lambda^{(1)}(h)= {\hat{\lambda}_1(h)}\, .\;\;\;\;\;\;
%\eea

\subsubsection{Transition in a twice longer interval}
\label{sec:two_step}
Consider now the interval $[0,2]$. We could have extended the method of Section~\ref{sec:one} and used the eigenvalue (square root of it) for the transition $s_0 \to s_2$ with transition density expressed in \eqref{shepp_form9} with $n=2$. This would improve Approximation 1 but this improvement is only marginal. Instead, we will use another approach: we consider the transition $s_1 \to s_2$ but use the interval $[0,1]$ just for  setting up the initial condition for observing $S(t)$ at  $t \in [1,2]$.

For $n=2$, the expression \eqref{shepp_form7} for the joint probability density function for the values $S(0),S(1), S(2)$ under the condition  $S(t)<h$ for all $t\in [0,2]$ has the form
\bea
p(s_0,s_1,s_{2})= \frac{1}{\varphi(s_0) F_2(h \big | s_0 )}
 %  \int_{-\infty}^h
   \det  \left(
                                                      \begin{array}{ccc}
                                                       \varphi(s_0)  & \varphi(s_0\!-\!h\!+\!s_1) &\varphi(s_0\!-\!2h\!+\!s_1\!+\!s_2) \\
                                                        \varphi(h) & \varphi(s_1)& \varphi(\!s_1\!+\!s_2\!-\!h) \\
                                                        \varphi(2h\!-\!s_1) & \varphi(h)& \varphi(s_2)
                                                      \end{array}
                                                    \right)
\, .
\eea

Denote by ${p}_1(z)$, $z<h,$  the `non-normalized' density of
 $S(1)$ under the condition $S(t)<h$ for all $t \in[0,1]$ that satisfies $\int_{-\infty}^{h}p_1(z)dz = F_1(h)$. Using \eqref{shepp_form3}, we obtain
 \bea
p_1(z) = \int_{-\infty}^{h}p^{(1 )}_h(x\to z )\varphi(x) dx = \Phi(h)\varphi(z) - \Phi(z)\varphi(h).
 \eea

Then the transition density from $x=s_1$ to $z=s_2$ under the condition  $S(t)<h$ for all $t\in [0,2]$ is achieved by integrating $s_0$ out
and renormalising the joint density:
\bea
{q}_h(x\to z )& =&
   \frac{1}{p_1(x)}\int_{-\infty}^h
   \det  \left(
                                                      \begin{array}{ccc}
                                                       \varphi(s_0)  & \varphi(s_0\!-\!h\!+\!x) &\varphi(s_0\!-\!2h\!+\!x\!+\!z) \\
                                                        \varphi(h) & \varphi(x)& \varphi(\!x\!+\!z\!-\!h) \\
                                                        \varphi(2h\!-\!x) & \varphi(h)& \varphi(z)
                                                      \end{array}
                                                    \right) ds_0
\\ &=&
\frac{1}{\Phi(h)\varphi(x) - \Phi(x)\varphi(h)} \det  \left(
                                                      \begin{array}{ccc}
                                                       \Phi(h)  & \Phi(x) &\Phi(\!x\!+\!z-h) \\
                                                        \varphi(h) & \varphi(x)& \varphi(\!x\!+\!z\!-\!h) \\
                                                        \varphi(2h\!-\!x) & \varphi(h)& \varphi(z)
                                                      \end{array}
                                                    \right) \, .
\eea
%
%
%\subsubsection{Two-step transition}
%\label{sec:two_step}
%
%
%In the case $n=2$ we obtain from \eqref{shepp_form2}:
% \bea
%p^{(2)}_h(x\to z ) = \frac{1}{\varphi(x)}  \int_{-\infty}^h \det  \left(
%                                                      \begin{array}{ccc}
%                                                       \varphi(x)  & \varphi(x\!-\!h\!+\!s_1) &\varphi(x\!-\!2h\!+\!s_1\!+\!z) \\
%                                                        \varphi(h) & \varphi(s_1)& \varphi(\!s_1\!+\!z\!-\!h) \\
%                                                        \varphi(2h\!-\!s_1) & \varphi(h)& \varphi(z)
%                                                      \end{array}
%                                                    \right) ds_1
%\eea
% with $z=s_2<h$. By integrating  $s_1$ out we obtain
%\bea \label{shepp_form5}
%p^{(2)}_h(x\to z ) = \frac{\varphi(h) \varphi(z)\Phi\left( \frac{v}{\sqrt{3}}\right)} {\sqrt{3}\varphi\left( \frac{v}{\sqrt{3}}\right) }
%\left[ 1\!-\!e^{hv-xz}\right]\!+\!\frac{ \varphi(h)[  \varphi(h)\Phi(v)\!-\!\Phi(x) \varphi(z)] }{\varphi(x)}\!+\!\Phi(h) \varphi(z)\!-\!\Phi(z) \varphi(h)\, ,
%\eea
%where $v=x+z-h$.
%

Let $\lambda_2(h)$ be the largest eigenvalue of
 the integral operator with kernel $q_h$:
 \bea
\lambda_2(h) q(z) = \int_{-\infty}^{h}q(x)q_h(x\to z ) dx, \text{  } z<h\, , \label{2.27}
\eea
where  eigenfunction $q(x)$ is some probability density on $(-\infty,h]$. Similarly to the case $n=1$,  $\lambda_2(h)$ is simple, real and positive eigenvalue
of the operator with kernel $K(x,z)=q_h(x\to z )$
 and  the eigenfunction $q(x)$ can be chosen as  a probability density.

In numerical examples below we approximate  $\lambda_2(h)$ using  the methodology described in \cite{Quadrature}, p.154. It is based on the Gauss-Legendre  discretization of  the interval $[-c,h]$, with some large $c>0$, into an $N$-point set $x_1, \ldots, x_N$ (the $x_i$'s are the roots of the \mbox{$N$-th} Legendre polynomial on $[-c,h]$), and the use of the Gauss-Legendre weights $w_i$ associated with points $x_i$; $\lambda_2(h)$ and $q(x)$ are  then approximated by the largest eigenvalue and associated eigenvector of the matrix
$
D^{1/2}AD^{1/2},
$
where $D = \text{diag}({w}_i)$, and  $A_{i,j} = q_h(x_i\to x_j )$. If $N$ is large enough then the resulting approximation $\hat{\lambda}_2(h)$ to $\lambda_2(h)$ is arbitrarily accurate.
\bea
\mbox{\noindent{\bf Approximation 2:}}\;F_T(h) \simeq  F_2(h) \left[ {\hat{\lambda}_2(h)} \, \right]^{T-2} \, \;\mbox{($T\geq 2$);} \;\;\Lambda^{(2)}(h) = - \log \hat{\lambda}_2(h) \, , \; \lambda^{(2)}(h)= {\hat{\lambda}_2(h)}\, .\;\;\;
\eea
%
%\noindent{\bf Approximation 2:}
%\bea
%F_T(h) \simeq  F_2(h) \left[ {\hat{\lambda}_2(h)} \, \right]^{T-2} \, \;\;\mbox{($T\geq 2$);}\;\; \;\;\Lambda^{(2)}(h) = - \log \hat{\lambda}_2(h) \, , \;\; \lambda^{(2)}(h)= {\hat{\lambda}_2(h)}\, .\;\;\;\;\;\;
%\eea

\subsubsection{Quality of Approximations 1 and 2}\label{Quality_of_l1}
Approximation 1 is more accurate than Approximation 0 but it is still not accurate enough. This is related to the fact that the process $S(t)$ is not Markovian and the behaviour of $S(t)$ on the interval $[i,i+1]$ depends on all  values of $S(t)$ in the interval $[i-1,i]$ and not only on the value $s_i=S(i)$, which is a simplification we used for derivation of Approximation~1. Approximation~2 corrects the bias of Approximation 1 by considering twice longer intervals $[i-1,i+1]$ and using the behaviour of $S(t)$ in the first half of the interval $[i-1,i+1]$ just for setting up the initial condition at $[i,i+1]$. As shown in Section~\ref{sim_section}, Approximation 2 is much more accurate than Approximations~0~and~1.
The approximations developed in the following section also carefully consider the dependence of $S(t)$ on its past; they  could be made arbitrarily accurate (on expense of increased computational complexity).

\subsection{Main approximations } \label{remaining_approx}

As mentioned above, the behaviour of $S(t)$ on the interval $[i,i+1]$ depends on all  values of $S(t)$ in the interval $[i-1,i]$ and not only on the value $s_i=S(i)$. The exact value of the Shepp's constant  $\lambda(h)$ can be defined as the limit (as $i \to \infty$) of the probability that $S(t)<h$ for all $t \in [i,i+1]$ under the condition
%that $S(t)$ has reached the stationary behaviour on the interval $[i-1,i]$ conditionally that
$S(t)<h$ for all $t \leq i$. Using the formula for conditional probability, we obtain
\be
\label{eq:lamb1}
\lambda(h) = \lim_{i\to \infty} F_i(h)/ F_{i-1}(h)\, .
\ee
Waiting  a long time without reaching $h$ is not numerically possible and is not what is really required for computation of  $\lambda(h)$. What we need is for the process $S(t)$ to (approximately) reach the stationary behaviour in the interval $[i-1,i]$ under the condition $S(t)<h$ for all $t<i$. Since the memory of $S(t)$ is  short (it follows from the representation $S(t)=W(t)-W(t+1)$, where $W(t)$ is the standard Wiener process), this stationary behaviour of $S(t)$ is practically achieved for very small~$i$, as is  seen from numerical results of  Section~\ref{sim_section}. Moreover, since ratios
 $F_i(h)/ F_{i-1}(h)$ are very close to  $F_i(h|x_h)/ F_{i-1}(h|x_h)$ for $i\geq 2$,  we can use ratios $F_i(h|x_h)/ F_{i-1}(h|x_h)$ in \eqref{eq:lamb1} instead. Here $x_h= - \varphi(h)/\Phi(h) $ is the mean of the truncated normal distribution with density $\varphi(x)/\Phi(h)$, $x\leq h$. For computing the approximations, it   makes  integration easier.
Note also another way of justifying the approximation $\lambda(h) \simeq F_i(h)/ F_{i-1}(h)$:  divide \eqref{eq:Shepp_constant1} with $T=i$ by \eqref{eq:Shepp_constant1} with $T=i-1$.

The above considerations give rise to several  approximations formulated below. We start with simpler approximations which are easy to compute and end up with  approximations which are extremely accurate but are harder to compute. Approximation 7 is very precise, see Table~\ref{Lambda_app}.
 However, we would not recommend extremely accurate Approximations 6 and 7 since Approximations 4 and~5 are already very accurate, see Tables 1 and 2, but are much easier to compute. Approximation~3, the simplest in the family, is also quite accurate. Note that all approximations for $F_T(h)$ can be applied for any $T\geq 2$.

\bea
\mbox{\noindent{\bf Approximation 3:}}\;\;F_T(h) &\simeq&  F_2(h) \left[ \lambda^{(3)}(h)  \right]^{T-2} \,, \;\;\mbox{where } \;\;\lambda^{(3)}(h) = {F_2(h|x_h)}/{ F_{1}(h|x_h)}  \, .\;\;\;\;\;\;
\\
\mbox{\noindent{\bf Approximation 4:}}\;\;F_T(h) &\simeq&  F_2(h) \left[ \lambda^{(4)}(h)  \right]^{T-2} \,, \;\;\mbox{where } \;\;\lambda^{(4)}(h) = {F_2(h)}/{ F_{1}(h)}   \, .\;\;\;\;\;\;
\\
\mbox{\noindent{\bf Approximation 5:}}\;\; F_T(h) &\simeq&  F_2(h) \left[ \lambda^{(5)}(h)  \right]^{T-2} \,, \;\;\mbox{where } \;\;\lambda^{(5)}(h) = {F_3(h|x_h)}/{ F_{2}(h|x_h)}   \, .\;\;\;\;\;\;
\\
\mbox{\noindent{\bf Approximation 6:}}\;\; F_T(h) &\simeq&  F_3(h) \left[ \lambda^{(6)}(h)  \right]^{T-3} \,, \;\;\mbox{where } \;\;\lambda^{(6)}(h) = {F_4(h|x_h)}/{ F_{3}(h|x_h)}  \, .\;\;\;\;\;\;
\\
\mbox{\noindent{\bf Approximation 7:}}\;\; F_T(h) &\simeq&  F_4(h) \left[ \lambda^{(7)}(h)  \right]^{T-4} \,, \;\;\mbox{where } \;\;\lambda^{(7)}(h) = {F_4(h)}/{ F_{3}(h)}  \, .\;\;\;\;\;\;
\eea

Numerical complexity of these approximation is related to the necessity of computing either $F_n(h|0)$ or $F_n(h)$  for suitable $n$. It follows from \eqref{shepp_form2} that $F_n(h|0)$  is an $n$-dimensional   integral. Consequently,  $F_n(h)$  is an $(n+1)$-dimensional   integral. In both cases, the dimensionality of the integral can be  reduced by one, respectively to $n-1$ and $n$, with no further analytical reduction possible.  In view of results of Sections~\ref{sec:app} and \ref{sec:app1}, computation of Approximations 3 and 4 is easy, computation of Approximation~5 requires numerical evaluation of a one-dimensional integral (which is not hard) but to compute Approximation~7 we need to approximate a three-dimensional integral, which has to be done with high precision as otherwise Approximation 7 is not worth using: indeed, Approximations~4--6 are almost as good but are much easier to compute. As  Approximation~7 provides us with the values which are  practically   indistinguishable from the true values of $\lambda(h)$, we use Approximation~7 only for
the assessment of the accuracy of other approximations and do not recommend using it in practice.

\subsection{Consistency of  approximations when $h$ is large}

Assume that $h \to \infty$. We shall show that Approximations 3-7 for $\Lambda(h)$ give consistent results with Approximation 0 which is
$\Lambda^{(0)}(h) = h\varphi(h)$.

Roughly, this consistency follows if we simply use $\Lambda^{(0)}(h)$ for $\Lambda(h)$ in \eqref{eq:Shepp_constant1} and then substitute the asymptotically correct values of $F_{i-1}(h)$ and $F_i(h)$ in
$\Lambda(h) \simeq \log F_{i-1}(h)- \log F_{i}(h)$. Similar argument works in the case  $\Lambda(h) \simeq \log F_{i-1}(h|x_h)- \log F_{i}(h|x_h)$.

Consider now Approximation 4 for $\Lambda(h)$, which is
$\Lambda^{(4)}(h) = \log F_1(h) - \log F_2(h)$.
From explicit formulas \eqref{eq:F1} and \eqref{eq:F2} for $F_1(h)$ and $ F_2(h)$ we obtain
\be
\label{eq:f1}
F_1(h)&=& 1- \left(h + \frac{2}{h} + O\left(\frac1{h^{3}}\right)\right) \varphi(h) \, , \;\;\;h \to \infty\, ,\\
F_2(h)&=& 1- \left(2h-4 - \frac{2}{h} + O\left(\frac1{h^{2}}\right)\right) \varphi(h) \, , \;\;\;h \to \infty\, ,
\label{eq:f2}\ee
Expansion \eqref{eq:f1} of \eqref{eq:F1} is straightforward. To obtain \eqref{eq:f2} from \eqref{eq:F2} we observe as $h \to \infty$:
$$
\Phi^3(h)= 1- \left(\frac{3}{h} + O\left(\frac1{h^{3}}\right) \right)\varphi(h)\, ; \;\;\; 2 h\varphi(h)\Phi^2(h)=
\left(2h - 4  + O\left(\frac1{h^{2}}\right) \right)\varphi(h)
$$
and
$$
\int_{-\infty}^{h}\Phi(y)^2\varphi(2h-y)dy=\left(\frac1{h} + O\left(\frac1{h^{2}}\right) \right)\varphi(h)\, ;
$$
all other terms in \eqref{eq:F2} converge to zero (as $h \to \infty$) faster than $\varphi(h)/h^2$. Using the expansion $\log(1-x) = -x + O(x^{-2})$ as $x \to 0$, this gives
\bea
\Lambda^{(4)}(h) = \log F_1(h) - \log F_2(h) = \left(h - 4 - \frac{4}{h}  + O\left(\frac1{h^{2}}\right) \right)\varphi(h)\;\;\mbox{as $h \to \infty$}\, .
\eea
This is fully consistent with approximation  $\Lambda^{(0)}(h)$  and all the discussion of Section~\ref{sec:21}. However, there is no guarantee that the constant 4 above is the correct constant in the asymptotic relation
\bea
\Lambda(h) =  \left(h - {\rm const}  + O\left(\frac1{h}\right) \right)\varphi(h)\;\;\mbox{as $h \to \infty$}\,
\eea
provided this asymptotic relation holds.

%
%
%%\noindent{\bf Approximation 3.}
%\bea
%\mbox{\noindent{\bf Approximation 3.}}\;\;F_T(h) \simeq  F_2(h) \left[ \lambda^{(3)}(h)  \right]^{T-2} \,, \;\;\mbox{where } \;\;\lambda^{(3)}(h) = {F_2(h|0)}/{ F_{1}(h|0)}  \, .\;\;\;\;\;\;
%\eea
%
%
%%\noindent{\bf Approximation 4.}
%\bea
%\mbox{\noindent{\bf Approximation 4.}}\;\;F_T(h) \simeq  F_2(h) \left[ \lambda^{(4)}(h)  \right]^{T-2} \,, \;\;\mbox{where } \;\;\lambda^{(4)}(h) = {F_2(h)}/{ F_{1}(h)}   \, .\;\;\;\;\;\;
%\eea
%
%
%%\noindent{\bf Approximation 5.}
%\bea
%\mbox{\noindent{\bf Approximation 5.}}\;\; F_T(h) \simeq  F_2(h) \left[ \lambda^{(5)}(h)  \right]^{T-2} \,, \;\;\mbox{where } \;\;\lambda^{(5)}(h) = {F_3(h|0)}/{ F_{2}(h|0)}   \, .\;\;\;\;\;\;
%\eea
%
%%\noindent{\bf Approximation 6.}
%\bea
%\mbox{\noindent{\bf Approximation 6.}}\;\; F_T(h) \simeq  F_3(h) \left[ \lambda^{(6)}(h)  \right]^{T-3} \,, \;\;\mbox{where } \;\;\lambda^{(6)}(h) = {F_3(h)}/{ F_{2}(h)}  \, .\;\;\;\;\;\;
%\eea
%
%
%%\noindent{\bf Approximation 7.}
%\bea
%\mbox{\noindent{\bf Approximation 7.}}\;\; F_T(h) \simeq  F_4(h) \left[ \lambda^{(7)}(h)  \right]^{T-4} \,, \;\;\mbox{where } \;\;\lambda^{(7)}(h) = {F_4(h)}/{ F_{3}(h)}  \, .\;\;\;\;\;\;
%\eea
%

%\section{Quality of approximations  and comparison of values of $\Lambda(h)$ for different processes}
\subsection{Numerical results}\label{sim_section}
\label{sec:3}
In this section we discuss the quality of approximations introduced in Section~\ref{sec:2}. In Table~\ref{lambda_approx}, we present the values of $\lambda^{(i)}(h), i=0,1,\ldots,7,$ for a number of different $h$; see also Table~\ref{Lambda_app} in Appendix.
As mentioned above,  $\lambda^{(7)}(h)$ is practically the true  $\lambda(h)$ and therefore we compare all other approximations against $\lambda^{(7)}(h)$. In Table~\ref{Relative_error_table} we present the relative errors of all other approximations against $\lambda^{(7)}(h)$; that is, the values $\lambda^{(i)}(h)/\lambda^{(7)}(h)-1$ for $i=0,1,\ldots,6$. From these two tables we see that Approximations 2-7 are very accurate. Moreover, we have made large-scale simulation studies  where we have estimated values of $F_T(h)$ for different $h$ using  $10^6$ trajectories of $S(t)$ and all  approximations for $F_T(h)$ considered above. Visually, Approximations 5-7 are virtually exact (the approximations are always well inside the confidence bounds computed from the simulations) for all $h\geq 0$ and also Approximations 2-4 are visually undistinguishable from them for $h\geq 0.5$. We do not provide corresponding plots as these plots are not informative.

\begin{table}[h]
\begin{footnotesize}
\begin{tabular}{|c||c|c|c|c|c|c|c|c|c|c|}
\hline
\!\! & $h$=0& $h$=0.5& $h$=1  & $h$=1.5& $h$=2& $h$=2.5& $h$=3& $h$=3.5& $h$=4\\
\hline
\!\!\!\!  $\lambda^{(0)}(h)$ \!\!\!\!  & 1.000000  & 0.838591  &  0.785079  & 0.823430  & 0.897644   & 0.957126  & 0.986792  &0.996950  &  0.999465   \\
\!\!\!\!  $\lambda^{(1)}(h)$ \!\!\!\!   & 0.250054  &  0.413754   &  0.596156 & 0.762590 & 0.885025   & 0.955674 &  0.986738  &0.996958  &0.999466   \\
\!\!\!\!  $\lambda^{(2)}(h)$ \!\!\!\!   & 0.201909& 0.366973  & 0.563246 & 0.746457  & 0.879719 &0.954522  & 0.986566  &0.996939    &0.999464    \\
\!\!\!\!  $\lambda^{(3)}(h)$ \!\!\!\!   &0.199421   &  0.366664  & 0.564851  & 0.747979   & 0.880220    &0.954529 & 0.986532 &  0.996930 & 0.999463   \\
\!\!\!\!  $\lambda^{(4)}(h)$ \!\!\!\!  &0.200045  & 0.365730 & 0.562888&0.746559 &0.879831 &0.954556 &0.986570  & 0.996939 &  0.999464\\
\!\!\!\!  $\lambda^{(5)}(h)$ \!\!\!\!  &0.202269  & 0.368099 & 0.564446 & 0.747143 & 0.879943 & 0.954564 &0.986571 & 0.996939  & 0.999464   \\
\!\!\!\!  $\lambda^{(6)}(h)$ \!\!\!\!   &0.202455  &0.368100  & 0.564377 & 0.747118 & 0.879945 &0.954566  &0.986571  &0.996939  & 0.999464    \\
\!\!\!\!  $\lambda^{(7)}(h)$ \!\!\!\!   &0.202434 & 0.368082 &0.564371   & {0.747118} &{0.879945} &{0.954566} & {0.986571} &{0.996939}&{0.999464}   \\
\hline
\end{tabular}
\caption{$\lambda^{(i)}(h), i=0,1,\ldots,7,$ for different $h$.}
\label{lambda_approx}
\end{footnotesize}
\end{table}

\begin{table}[h]
\begin{footnotesize}
\begin{tabular}{|c||c|c|c|c|c|c|c|c|c|c|}
\hline
\!\!  &$h=0$& $h$=0.5& $h$=1  & $h$=1.5& $h$=2& $h$=2.5& $h$=3& $h$=3.5& $h$=4\\
\hline
\!\!\!\!  $\lambda^{(0)}(h)$ \!\!\!\!   & 3.94e+00 &  1.28e+00 &3.91e-01  &1.02e-01  & 2.01e-02   &2.68e-03  &  2.25e-04  &   1.16e-05     &  6.12e-07    \\
\!\!\!\!  $\lambda^{(1)}(h)$ \!\!\!\!     &2.35e-01  &1.24e-01  &5.63e-02   &2.07e-02  &  5.77e-03  & 1.16e-03 &1.69e-04    & 1.93e-05       & 1.88e-06     \\
\!\!\!\!  $\lambda^{(2)}(h)$ \!\!\!\!    & -2.59e-03 & -3.01e-03  & -1.99e-03  & -8.84e-04 & -2.57e-04   & -4.56e-05 &  -4.61e-06  & -2.56e-07       & -7.82e-09     \\
\!\!\!\!  $\lambda^{(3)}(h)$ \!\!\!\!    &-1.49e-02  & -3.85e-03  & 8.51e-04  & 1.15e-03 & 3.12e-04   & -3.84e-05 & -3.88e-05   & -9.36e-06       &-1.28e-06      \\
\!\!\!\!  $\lambda^{(4)}(h)$ \!\!\!\!    &-1.18e-02  &-6.39e-03   &-2.63e-03   &  -7.48e-04&  -1.29e-04  & -1.09e-05 & -2.06e-07 &2.27e-08        &1.35e-09      \\
\!\!\!\!  $\lambda^{(5)}(h)$ \!\!\!\!    &-8.13e-04  & 4.71e-05  &1.33e-04  &3.32e-05  & -2.49e-06   & -1.57e-06& -1.34e-07   & -2.20e-09       & 9.09e-11     \\
\!\!\!\!  $\lambda^{(6)}(h)$ \!\!\!\!  & 1.03e-04 &5.02e-05  & 1.09e-05  &  1.88e-07 &-1.83e-07    & -3.22e-11 &   4.12e-11  &2.86e-11        &   6.09e-12   \\
\hline
\end{tabular}
\end{footnotesize}
\caption{Relative errors of $\lambda^{(i)}(h), i=0,1,\ldots,6,$ against $\lambda^{(7)}(h)$.}
\label{Relative_error_table}
\end{table}

\begin{figure}[!h]
%\label{table1}
\label{fig:lambda12}
\begin{center}
\begin{subfigure}{0.48\linewidth} \centering
     \includegraphics[scale=0.2]{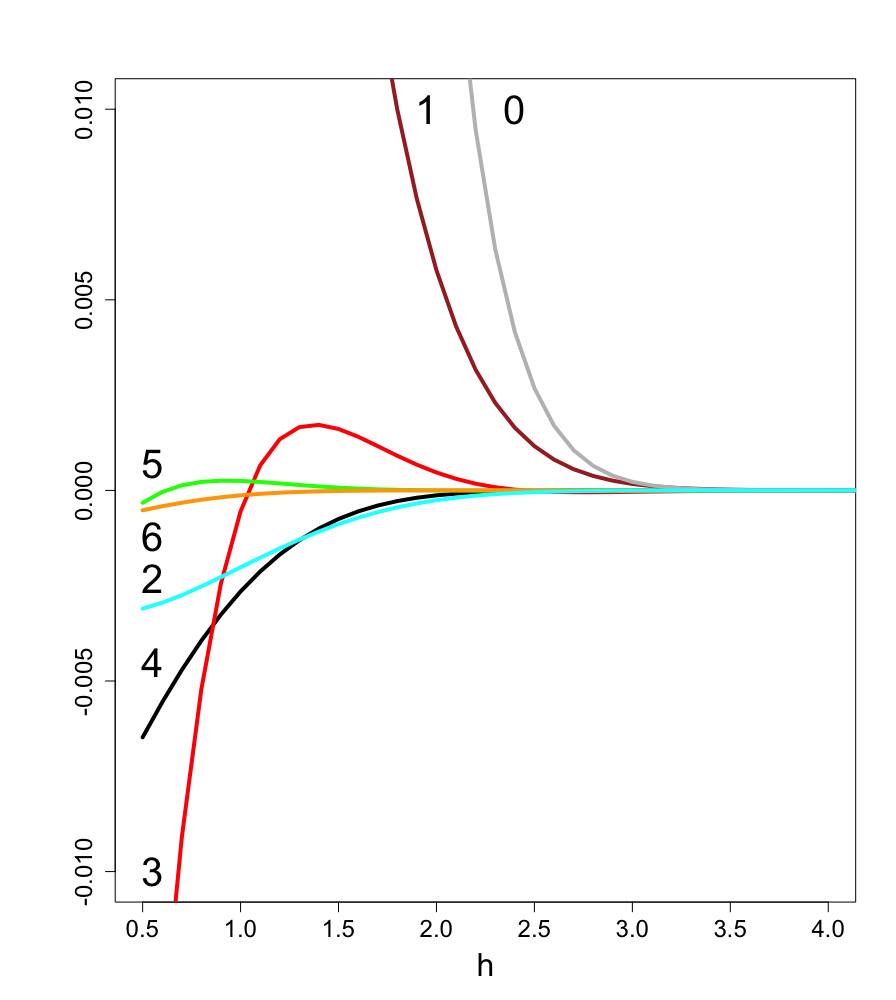}
     \caption{Relative errors of $\lambda^{(i)}(h)$, $i=0,\ldots,6$, against $\lambda^{(7)}(h)$  } \label{fig:RE}
\end{subfigure}
\hspace{0.3cm}
 \begin{subfigure}{0.48\linewidth} \centering
     \includegraphics[scale=0.2]{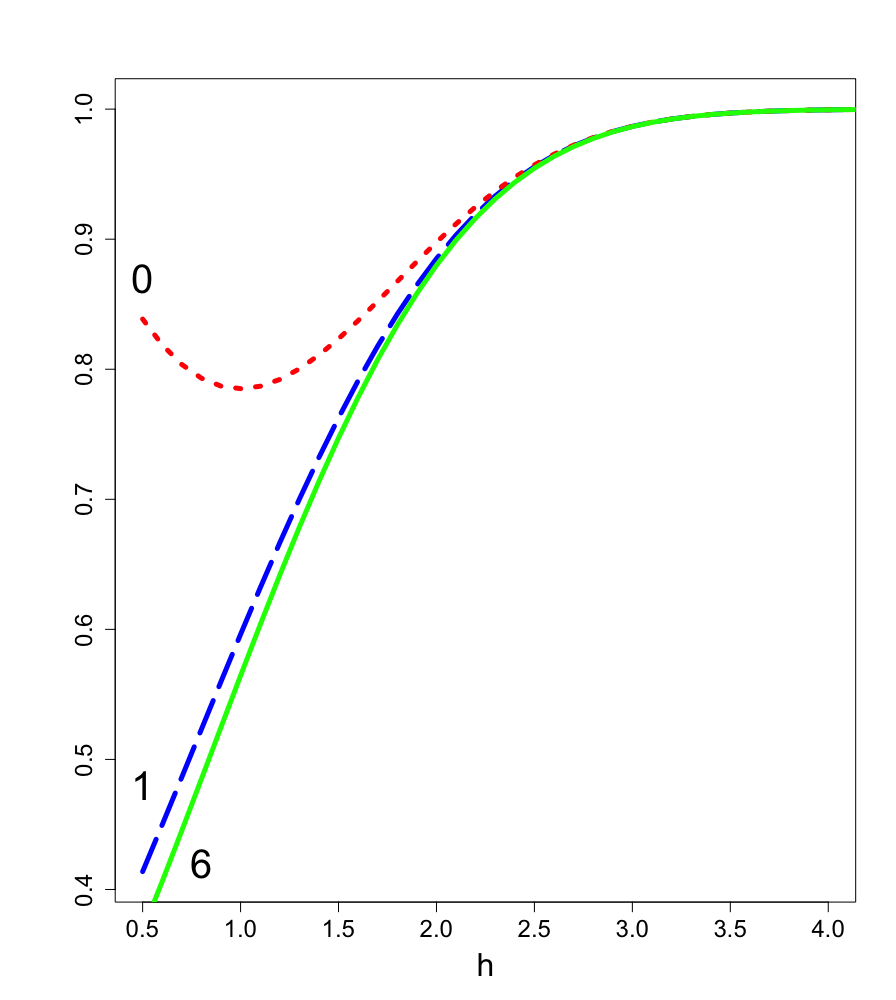}
     \caption{$\lambda^{(0)}(h)$ (dotted red), $\lambda^{(1)}(h)$ (dashed blue)  and  $\lambda^{(6)}(h)$ (solid green) } \label{fig:lambda}
\end{subfigure}
    \end{center}
    \caption{Approximations and their relative errors as functions of $h$.}
\end{figure}

A plot of the relative errors can be seen in Figure~\ref{fig:RE}, where the number next to the line corresponds to the approximation.
 Approximations 2,4 and 7  suggest very accurate lower bounds for the true $\lambda(h)$. Approximations 0 and 1 appear to provide upper bounds for $\lambda(h)$ for all $h$.

As mentioned in Section~\ref{Quality_of_l1}, Approximation 1 is not as accurate as Approximations 2--7 because it does not adequately take into account  the non-Markovianity of  $S(t)$. In Figure~\ref{fig:lambda} we have plotted  $\lambda^{(0)}(h)$ (dotted red line), $\lambda^{(1)}(h)$ (dashed red line)  and  $\lambda^{(6)}(h)$ (solid green line) for a range of interesting $h$. Visually, all $\lambda^{(i)}(h)$ with $i=2,4,5,6,7$
would be visually  indistinguishable from each other on the plot in Figure~\ref{fig:lambda} and $\lambda^{(3)}(h)$ would be very close to them.
 The number next to the line corresponds to which approximation was used.

%\AZ{ WILL DO SOMETHING WITH THIS TEXT LATER.
%In Figure~\ref{sim_lam} we illustrate the rate of convergence of $   [\log F_n(h)]/n $ in \eqref{eq:Shepp_constant} to  the %Shepp's constant $\Lambda(h)$ as $n$ increases. The dotted black lines correspond to simulation results obtained by computing the probability $F_n(h)$ with 100,000 simulations for $h=1.5,2$ and $3$. The solid red lines correspond to Approximation~4 for $\Lambda(h)$ with the chosen $h$. Figure~\ref{sim_lam} demonstrates how accurate this computationally cheap approximation is and also demonstrates the importance of the multiplying constant $F_2(h)$ in Approximation~4 to correct for the non-linear behaviour seen for small~$n$. In Figure~\ref{shepp_con} we investigate the rate of convergence and accuracy of convergence to $\Lambda(h)$ using all approximations, where we have fixed $h=3$. In this figure, Approximations 4, 5, 6 and 7 produce results that are visually  indistinguishable  to Approximation 2 and hence are not plotted. Therefore, in this figure,  Approximation~2 can be considered as giving the true Shepp's constant $\Lambda(h)$.
% }

\section{Comparison of the upper tail asymptotics for the Slepian process against some other stationary Gaussian processes}
\label{sec:4}
Consider the following three stationary Gaussian processes.
\begin{enumerate}
 \item $\xi_1(t)$  ($ t \geq 0$) is the Ornstein-Uhlenbeck process with mean 0, variance 1 and correlation function $\rho_1(t)=\exp(-|t|)$.
 \item Let $a>0$ be fixed real number and set $\alpha = (1+a+a^2)/(2+2a+a^2)$. Then, if $W(t)$ denotes the standard Wiener process, we define the process $\xi_2(t)$ ($ t \geq 0$) as follows:
\bea
\xi_2(t) &=& \frac{1}{\sqrt{1+a+a^2}} \left\{ (1+a)W(t+2\alpha) - aW\left(t+{\alpha}\right) - W(t)  \right\}.\nonumber
\eea
The process $\xi_2(t)$ has mean 0, variance 1 and correlation function
\begin{equation*}
 \rho_2(t)=
  \begin{cases}
    1-|t|, & \text{for } 0 \leq |t| \leq {\alpha} \\
    \frac{(1+a)(2\alpha-|t|)}{1+a+a^2}, & \text{for }  {\alpha} \leq |t| \leq 2\alpha\\
    0& \text{for }   |t| \ge 2\alpha\, .
  \end{cases}
 \end{equation*}
 \item  Let $c\ge1$ be a fixed real number and set $\beta=1/(c+2)$. Define the process  $\xi_3(t)$ by
 \begin{eqnarray*}
\xi_3(t) &=& \frac{{1}}{\sqrt{1+c^2}} \left\{ W(t+1) + cW\left(t+(c+1)\beta\right) - cW\left(t+\beta \right) -W(t) \right\}\nonumber.
\end{eqnarray*}
The process $\xi_3(t)$ has  mean 0, variance 1 and correlation function
\begin{equation*}
\rho_3(t)=
  \begin{cases}
    1-|t|, & \text{for } 0 \leq |t| \leq \beta \\
    \frac{(1+c)(1+c^2\beta -|t|(1+c))}{1+c^2}& \text{for } \beta \leq |t| \leq c\beta \\
     \frac{1+c+c^2\beta-|t|(1+2c)}{1+c^2} & \text{for }  c\beta \leq |t| \leq (c+1)\beta,\\
        \frac{1-|t|}{1+c^2}, & \text{for }    (c+1)\beta\leq |t| \leq 1\\
    0, & \text{for }     |t| \geq 1.
  \end{cases}
 \end{equation*}

%
%\begin{equation*}
%\rho_3(t)=
%  \begin{cases}
%    1-|t|, & \text{for } 0 \leq |t| \leq \frac{1}{3} \\
%    \frac{7}{6}-\frac{3}{2}|t| & \text{for }  \frac{1}{3} \leq |t| \leq \frac{2}{3} \\
%    \frac{1}{2} - \frac{|t|}{2}& \text{for }  \frac{2}{3}  \leq |t| \leq 1,\\
%    0, & \text{for }     |t| \geq 1.
%  \end{cases}
% \end{equation*}
\end{enumerate}

\begin{figure}[!h]
\label{table1}
\begin{center}
\begin{subfigure}{0.48\linewidth} \centering
     \includegraphics[scale=0.22]{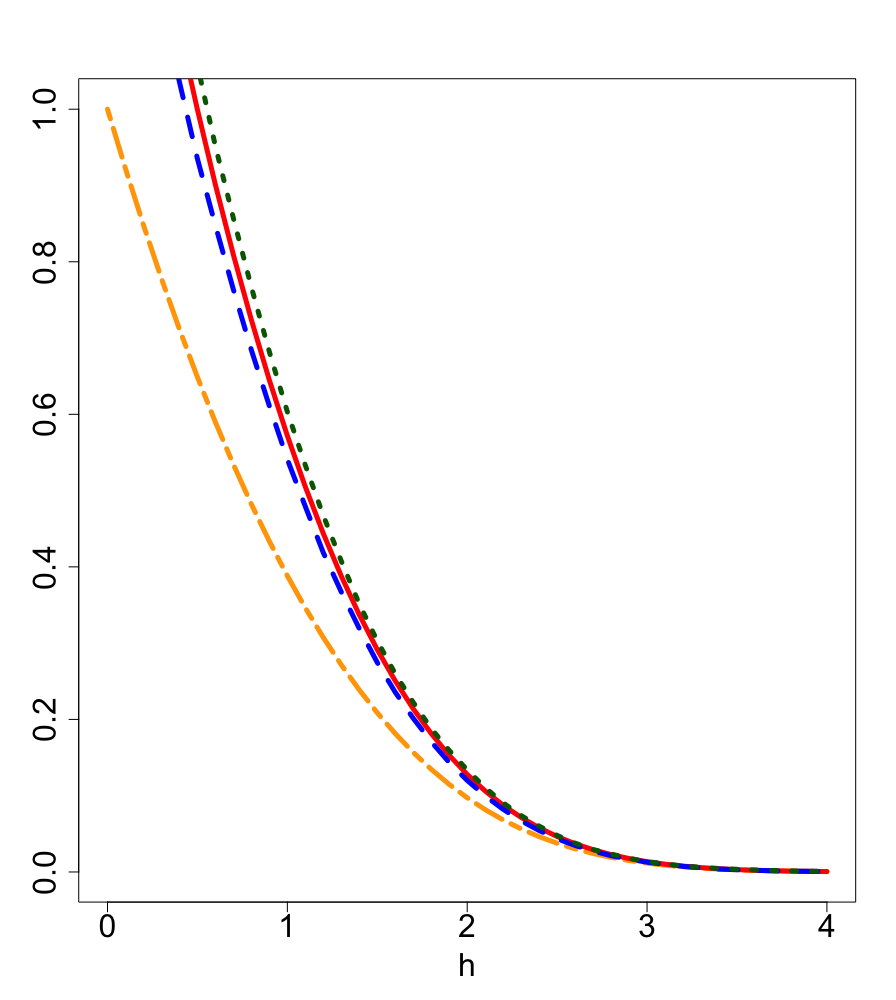}
     \caption{ Shepp's constant $\Lambda(h)$ and $\Lambda_i(h)$, $i=1,2,3$ } \label{fig:Lambda_comparison}
\end{subfigure}
\hspace{0.3cm}
 \begin{subfigure}{0.48\linewidth} \centering
     \includegraphics[scale=0.22]{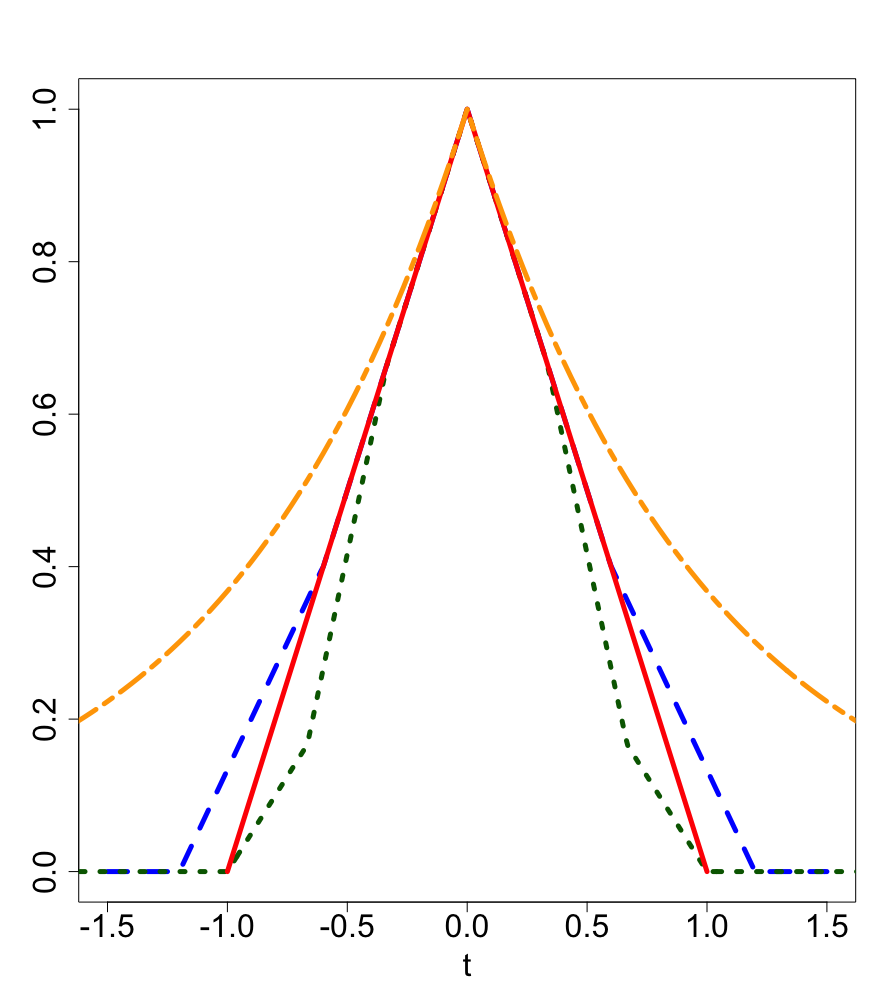}
     \caption{Correlation functions $\rho(t)$ and $\rho_i(t)$, $i=1,2,3$} \label{fig:correlation_comparison}
\end{subfigure}
    \end{center}
    \vspace{-0.4cm}
    \caption{Comparison of the upper tail asymptotics for several Gaussian stationary  process  }
\end{figure}

It follows from  \cite[Theorem 3]{Mehr}, the above three processes provide a very good representation of the entire class of conditionally Markov stationary Gaussian processes. Indeed, there is only one  process in this class where $\alpha \neq 1$ in \eqref{eq:cov_f} (this is the process with covariance function $\rho(t)=\cos \omega t $ with $\omega \neq 0$) and the three types of processes we consider cover well the case where $\alpha=1$ and $C=1$  in \eqref{eq:cov_f} (the case  $C \neq 1$ reduces to the case $C=1$ by substituting  $h/C$ for $h$). For a graphical representation of the chosen covariance functions, see
Figure~\ref{fig:correlation_comparison}.

Below we compare Shepp's constant $\Lambda(h)$ defined in \eqref{eq:Shepp_constant} to similar quantities of the  processes $\{ \xi_i(t)\}$, ($i=1,2,3$)  defined above. More precisely, let
\bea\label{eq:bcp}
F_{T,i}(h ) := {\rm Pr}\left\{\max_{t \in [0,T]} \xi_i(t) < h\right\}\,,\;\; i=1,2,3.
\eea
We are interested in comparing Shepp's constant $\Lambda(h)$ with
\be\label{eq:general_shepp}
\Lambda_i(h) = -\lim_{T \to \infty} \frac1T   \log F_{T,i}(h),
\ee
for $i=1,2,3.$ Importantly, each process has $\mathbb{E}\xi_i(t)=0$, $\mathbb{E}\xi_i^2(t)=1$ and correlation function $\rho_i(t)=\mathbb{E}\xi_i(0)\xi_i(t)$ which satisfies $\rho_i'(0^+)=\frac{d}{dt}\rho_i(t)|_{t=0^+} =-1$, for $i=1,2,3$.

%\subsection{Existence and approximations for $\Lambda_i(h)$}
 The existence and evaluation of the constant $\Lambda_1(h)$ defined in \eqref{eq:general_shepp} for the Ornstein-Uhlenbeck process has been considered in \cite{beekman1975asymptotic}, where  it was shown that  $0<\Lambda_1(h)<1$ for all $h>0$, and that $\Lambda_1(h)$ is the root of a parabolic cylinder function (defined in \cite{beekman1975asymptotic}) closest to zero. It is also shown that $\lim_{h\to 0^{+}}\Lambda_1(h)=1$ and $\lim_{h\to \infty}\Lambda_1(h)=h\varphi(h)$. The existence of the constants  $\Lambda_2(h)$ and $\Lambda_3(h)$ follows from similar arguments for the existence of Shepp's constant $\Lambda(h)$. Moreover, the constants are approximated by the same methodology as Shepp's constant, namely:
 \be\label{eq:lambda_general}
 \lambda_i(h) = \exp(-\Lambda_i(h)) = \lim_{j\to \infty} F_{j,i}(h)/ F_{j-1,i}(h)\, .
 \ee
 The justification why we expect $ F_{j,i}(h)/ F_{j-1,i}(h)$ (with, say, $j\geq 3$) to be a good approximation of $\lambda_i(h)$ is related to the property of `fast loss of memory', which  processes $\xi_2(t)$ and $\xi_3(t)$ possess, as the process $S(t)$ does.
 In view of the complex structure of $\xi_2(t)$ and $\xi_3(t)$, the values of  $F_{T,2}(h )$ and $F_{T,3}(h )$  are evaluated via Monte Carlo simulations. In Figure~\ref{fig:Lambda_comparison}, we compare $\Lambda_i(h)$ with Shepp's constant $\Lambda(h)$  (red solid line). $\Lambda_1(h)$ (orange dot-dash line) has been computed as in \cite{beekman1975asymptotic} . $\Lambda_2(h)$ (blue dashed line) and $\Lambda_3(h)$ (dark green dotted line) have been approximated using \eqref{eq:lambda_general} with $j=3$. For $\Lambda_2(h)$ we have taken $a=1$ in the definition of $\xi_2(t)$ and for $\Lambda_3(h)$ we have taken $c=1$ in the definition of $\xi_3(t)$.  In Figure~\ref{fig:correlation_comparison}, we plot the correlation functions: $\rho(t)$ (red solid line); $\rho_1(t)$ (orange dot-dash line); $\rho_2(t)$ (blue dashed line); $\rho_3(t)$ (dark green dotted line). Note that the results obtained are fully consistent with the celebrated `Slepian's lemma',  a Gaussian comparison inequality, see Lemma 1 in \cite{slepian1962}. In our terms, Slepian's lemma says that if for two stationary Gaussian  processes with non-negative covariance functions $\rho_1$ and $\rho_2$ we have $\rho_1(t) \geq \rho_2(t)$ for all $t\geq 0$, then for the corresponding values of $\Lambda(h)$ we have $\Lambda_1(h)$ {$\leq \Lambda_2(h)$}, for all $h$. \\

%The inequality \eqref{eq:Slep_ineq} follows from the infamous `Slepian's lemma', which  is a Gaussian comparison inequality, see Lemma 1 in \cite{slepian1962}. It states that for Gaussian random variables $X = (X_1,\dots,X_n)$ and $Y = (Y_1,\dots,Y_n)$ in $\mathbb{R}^n$ satisfying $\mathbb{E}[X] = \mathbb{E}[Y] = 0$,
%$\mathbb{E}[X_i^2]=\mathbb{E}[Y_i^2]$, $i=1,\dots,n, \text{ and } \ \mathbb{E}[X_iX_j] \le \mathbb{E}[Y_i Y_j] $ for $i \neq j$,
%the following inequality holds for all real numbers $u_1,...,u_n$:
%$P[X_1 \le u_1, \dots, X_n \le u_n] \le P[Y_1 \le u_1, \dots, Y_n \le u_n ] $. To derive \eqref{eq:Slep_ineq} from the Slepian's lemma we can follow the steps of the proof (of a  particular case of \eqref{eq:Slep_ineq} where $h=0$) of Section 2.4 in \cite{slepian1962}; see also the Wikipedia article `Slepian's lemma'.
%
%Now we  use other Slepian's arguments; this time  from \cite[p.469]{slepian1962}

%\begin{figure}[!h]
%\label{table2}
%\begin{center}
% \begin{subfigure}{0.48\linewidth} \centering
%     \includegraphics[scale=0.2]{Sim_of_shepp_cont.png}
%     \caption{Rate and accuracy of convergence to  the Shepp's constant  $\Lambda(h)$ using simulations (dotted black line) and using Approximation~4 (solid red) for $h=1.5, 2$ and 3.}\label{sim_lam}
%   \end{subfigure}
%   \hspace{0.3cm}
%   \begin{subfigure}{0.48\linewidth} \centering
%     \includegraphics[scale=0.2]{Shepp_approximations.png}
%     \caption{Rate and accuracy of convergence to $\Lambda(h)$ with $h=3$ using all approximations (numbers next to line correspond to which approximation has been used).}\label{shepp_con}
%   \end{subfigure}
%\end{center}
%\caption{ }
%\end{figure}

\nopagebreak
\section{Appendix}

\subsection{Approximations for $\Lambda(h)$}
In Table~\ref{Lambda_app} we use Approximation~7 (our most accurate approximation) to approximate $\Lambda(h)$ over  increments 0.1 for $h$.
 Bold font indicates the  decimal places which we claim  accurate. Note that $h=0$ has been treated as a special case, see for example \cite{slepian1962} and
\cite{pitman2015slepian}. For $h=0$, instead of Approximation~7, we have  used the approximation $\Lambda^{(8)}(h) = -\log(F_5(h)/F_4(h))$; we do not recommend using this approximation in general because of its high complexity.\\

\begin{table}[H]
\begin{tabular}{|c|c||c|c||c|c||c|c||c|c|}
  \hline
  $h$ & $\Lambda(h) $ &$h$ & $\Lambda(h) $ &$h$ & $\Lambda(h) $ &$h$ & $\Lambda(h) $ &$h$ & $\Lambda(h) $  \\ \hline
  0.0 & \textbf{1.5972} & 0.8 & \textbf{0.724}0 & 1.6 & \textbf{0.2505}19 & 2.4 & \textbf{0.0578944} & 3.2 & \textbf{0.0077016} \\
  0.1 & \textbf{1.463}2 & 0.9 & \textbf{0.645}0 & 1.7 & \textbf{0.2139}29 & 2.5 & \textbf{0.0464986} & 3.3 & \textbf{0.0057244} \\
  0.2 & \textbf{1.336}5 & 1.0 & \textbf{0.572}0 & 1.8 & \textbf{0.18148}4 & 2.6 & \textbf{0.0370122} & 3.4 & \textbf{0.0042111}\\
  0.3 & \textbf{1.217}0 & 1.1 & \textbf{0.5051}& 1.9 & \textbf{0.15290}2 & 2.7 & \textbf{0.0291909} & 3.5 & \textbf{0.0030658} \\
  0.4 & \textbf{1.104}7 & 1.2 & \textbf{0.4438} & 2.0 & \textbf{0.127896} & 2.8 & \textbf{0.0228058} & 3.6 & \textbf{0.0022087} \\
  0.5 & \textbf{0.999}5& 1.3 & \textbf{0.3879} & 2.1 & \textbf{0.106178} & 2.9 & \textbf{0.0176462} & 3.7 & \textbf{ 0.0015747 } \\
  0.6 & \textbf{0.901}0& 1.4 & \textbf{0.3372} & 2.2 & \textbf{0.087460} & 3.0 & \textbf{0.0135203} & 3.8 & \textbf{0.0011109} \\
  0.7 & \textbf{0.809}2 & 1.5 & \textbf{0.2915} & 2.3 & \textbf{0.071458} & 3.1 & \textbf{0.0102561} & 3.9 & \textbf{  0.0007755 } \\
  \hline
\end{tabular}
\caption{Approximations for $\Lambda(h)$ with accurate decimal digits in bold.}
\label{Lambda_app}
\end{table}

\label{sec:app}
\subsection{An approximation for $F_2(h)$}
\label{sec:app1}

Using approximations for $\Phi(t)$, it is possible to approximate $F_2(h)$ very accurately. For example, using the approximation
(see  \cite{lin1989approximating})
\begin{eqnarray} \label{Lin_form}
\Phi(t) = \left\{ \begin{array}{ll}
    0.5\exp(0.717t - 0.416t^2) & \text{for } t\leq0\\
        1-0.5\exp(-0.717t - 0.416t^2) & \text{for } t>0\, ,\\
       \end{array} \right.
\end{eqnarray}
we obtain
\begin{eqnarray}\label{Corrected_Diffusion_explicit2}
F_2(h) &\cong&  \Phi(h)^3+\varphi({h})^2\Phi(h) +\frac{\varphi({h})^2}{2} \left[ (h^2-1)\Phi(h) + h\varphi(h) \right] - 2\varphi({h})\Phi({h})\left[ h\Phi(h)+\varphi(h) \right]   \nonumber \\
&+& \Phi(2{h})-\Phi({h})  -  \frac{0.5}{\sqrt{2\pi}}e^{-2{h}^2 } \bigg[  2J(0.916,b,{h}) - \frac{1}{2}J(1.332,b_1,{h})-\frac{1}{\sqrt{2\pi}}V(1.416,b,{h}) \nonumber \\
&+&  \frac{2}{\sqrt{2\pi}}V(1,b_2,{h}) + \frac{1}{\pi}K(1.5,b_2,{h})-\frac{1}{2} \left\{  K(1.332,b_3,0) -\frac{2}{\sqrt{2\pi}}U(1.416,b_4,0) \right \}  \bigg]\, ,\;\;\;\;\;\;\;\;\;\;
\end{eqnarray}
where
$ b = 2{h} - 0.717, b_1 = b-0.717, b_2=2{h},
b_3 = b+2.151, b_4 = b+1.434$,
\begin{eqnarray}
K(x,y,z) =  \frac{\sqrt{\pi}e^{y^2/(4x)}}{\sqrt{x}}\Phi \left(\frac{2xz-y}{\sqrt{2x}}  \right)\,, \;\; U(x,y,z)=\frac{1}{2x} \left[y K(x,y,z)- e^{z(y-xz)}\right]\, ,  \label{K_form}
%U(x,y,z) &=& \frac{\sqrt{\pi}ye^{y^2/(4x)}}{2x^{3/2}}\Phi \left(\frac{2xz-y}{\sqrt{2x}}   \right) - \frac{e^{z(y-xz)}}{2x}, \nonumber\\
\end{eqnarray}
%$U(x,y,z)=K(x,y,z)\cdot y/(2x) - \exp(z(y-xz))/2x$,
$J(x,y,z) = K(x,y,z) - K(x,y,0)$ and $ V(x,y,z) = U(x,y,z) - U(x,y,0).$\\\\
Table~\ref{F2_approx} shows that approximation \eqref{Corrected_Diffusion_explicit2} is very accurate across all $h$ of interest.
\begin{table}[h]
\begin{small}
\begin{tabular}{|c||c|c|c|c|c|c|c|c|c|c|}
\hline
\!\! & $h=0$&$h$=0.5& $h$=1  & $h$=1.5& $h$=2& $h$=2.5& $h$=3& $h$=3.5& $h$=4\\
\hline
\!\!\!\!  $F_2(h) $ \!\!\!\!   &0.018173 & 0.085014  & 0.250896  & 0.502268  &  0.744845    &0.900875   & 0.970790  & 0.993430  & 0.998866    \\
\!\!\!\!  \eqref{Corrected_Diffusion_explicit2} \!\!\!\!   &0.019548  &0.084687  &0.250203 &0.502097  &0.744837 &0.900875 &0.970790   &0.993430 &0.998866\\
\hline
\end{tabular}
\caption{Accuracy of approximation \eqref{Corrected_Diffusion_explicit2} for $F_2(h)$}
\label{F2_approx}
\end{small}
\end{table}

\subsection{Simplified form of $F_2(h|x_h)$ and its approximation }
Using \eqref{shepp_form}, for any $x_0\leq 0$,  we can express $F_2(h|x_0)$ as follows:
\bea
F_2(h|x_0) &=& \Phi(h)^2 +\frac{1}{\varphi(x_0)}\varphi(h)^2x_0\Phi(x_0)-\frac{1}{\varphi(x_0)}\varphi(h)\Phi(h)\Phi(x_0)-h\varphi(h)\Phi(h) \\
&+& \frac{1}{\varphi(x_0)}\int_{h}^{\infty}\varphi(y)\Phi(2h-y)\varphi(h+x_0-y) dy
 - \frac{1}{\varphi(x_0)}\int_{h}^{\infty}\Phi(h+x_0-y)\varphi(2h-y)\varphi(y) dy.
\eea

Using \eqref{Lin_form}, we  obtain the approximation $F_2(h|x_0) \cong \hat{F}$ where

\bea
\!\!\!\!\!\hat{F}\!\!\!\!\! &=& \Phi(h)^2 +\frac{1}{\varphi(x_0)}\varphi(h)^2x_0\Phi(x_0)-\frac{1}{\varphi(x_0)}\varphi(h)\Phi(h)\Phi(x_0)-h\varphi(h)\Phi(h) \nonumber\\
&+&\frac{1}{\sqrt{2}\varphi(x_0)} \varphi\left (\frac{h+x_0}{\sqrt{2}} \right)\left[ \Phi\left (2h\sqrt{2} - \frac{h+x_0}{\sqrt{2}} \right ) - \Phi\left (h\sqrt{2} - \frac{h+x_0}{\sqrt{2}} \right )  \right] \nonumber \\
&-& \frac{\varphi(h\!+\!x_0)}{2\sqrt{2\pi}\varphi(x_0)}e^{-1.664 h^2}  \{  e^{-1.434h} \left[ K(1.416,2.664h\!+\!x_0\!+\!0.717,2h)\!-\!K(1.416,2.664h\!+\!x_0\!+\!0.717,h)    \right] \;\;\;\;\; \nonumber\\
&-& e^{1.434h}[ K(1.416,2.664h+x_0\!-\!0.717,\infty ) - K(1.416,2.664h\!+\!x_0\!-\!0.717,2h)]   \}   \nonumber \\
&-& \frac{e^{0.717(h+x_0)-2h^2-0.416(h+x_0)^2}}{2\varphi(x_0){2\pi}} \bigg[ K(1.416,2h+0.832(h+x_0)-0.717,\infty)\nonumber \\
&-&K(1.416,2h+0.832(h+x_0)-0.717,h) \bigg ],\label{f_2_zero}\;\;\;\;\;\;\;\;\;\ \nonumber
\\\nonumber
\eea
where $K(x,y,z)$ is defined  in \eqref{K_form}.
Table~\ref{F2_z_approx} shows that approximation \eqref{f_2_zero} is very accurate for $x_0=x_h= -\varphi(h)/\Phi(h)$, for any  $h\geq 0$.
\begin{table}[h]
\begin{footnotesize}
\begin{tabular}{|c||c|c|c|c|c|c|c|c|c|c|}
\hline
\!\! &$h$=0& $h$=0.5& $h$=1  & $h$=1.5& $h$=2& $h$=2.5& $h$=3& $h$=3.5& $h$=4\\
\hline
\!\!\!\!  $F_2(h|x_0) $ \!\!\!\!   &  0.041459  &0.141066 & 0.337112 & 0.588949   & 0.803170   &  0.927924 & 0.979740  &  0.995608 &  0.999264  \\
\!\!\!\!  $\hat{F}$ \!\!\!\!   & 0.041942 & 0.139821   &  0.336115& 0.588695  & 0.803139  & 0.927922  & 0.979740   &  0.995608  &0.999264  \\
\hline
\end{tabular}
\caption{Accuracy of approximation $\hat{F}\simeq F_2(h|0)$}
\label{F2_z_approx}
\end{footnotesize}
\end{table}

\section{Conclusions}
\label{sec:conc}

In his seminal paper \cite{Shepp71}, L. Shepp derived explicit formulas
for $F_T(h ) = {\rm Pr}\left\{\max_{t \in [0,T]} S(t) < h   \right\}
$, the distribution of maximum of the so-called Slepian process $S(t)$.
As these explicit formulas are complicated, in the same paper L. Shepp
has introduced a constant  $\Lambda(h) = -\lim_{T \to \infty} \frac1T   \log F_T(h)$ (which we call Shepp's constant) measuring the rate of decrease of $F_T(h)$ as $T$ grows; L. Shepp  also  raised the question of constructing accurate approximations and bounds for this constant.
Until now, this question has not been adequately addressed. To answer it, we have constructed different  approximations for $F_T(h)$ (and hence for
$\Lambda(h)$). We have  shown in Section~\ref{sec:3} that at least some of these approximations are extremely accurate for all $h\geq 0$. We have also provided other approximations that are almost as good but are much simpler to compute.

\section*{Acknowledgement}
The authors are grateful to the AE and both reviewers whose comments much helped for improving the presentation.

%\cite{harper2017pickands}
%\cite{li2004lower}
%\cite{molchan2012survival}
%\cite{pitman2015slepian}

\bibliographystyle{unsrt}

%\cite{beekman1975asymptotic}
%\cite{Finch2004ornstein}
%\bibliographystyle{natbib}

\bibliography{changepoint}

\begin{thebibliography}{10}

\bibitem{Shepp71}
LA~Shepp.
\newblock First passage time for a particular {G}aussian process.
\newblock {\em The Annals of Mathematical Statistics}, pages 946--951, 1971.

\bibitem{landau1970supremum}
HJ~Landau and LA~Shepp.
\newblock On the supremum of a {G}aussian process.
\newblock {\em Sankhy{\=a}: The Indian Journal of Statistics, Series A}, pages
  369--378, 1970.

\bibitem{marcus1972sample}
M~Marcus and LA~Shepp.
\newblock Sample behavior of {G}aussian processes.
\newblock In {\em Proc. of the Sixth Berkeley Symposium on Math. Statist. and
  Prob}, volume~2, pages 423--421, 1972.

\bibitem{adler2007random}
RJ~Adler and J~Taylor.
\newblock {\em Random Fields and Geometry}.
\newblock Springer, 2007.

\bibitem{Aldous}
D~Aldous.
\newblock {\em Probability {A}pproximations via the {P}oisson {C}lumping
  {H}euristic}.
\newblock Springer Science \& Business Media, 1989.

\bibitem{Mehr}
CB~Mehr and JA~McFadden.
\newblock Certain properties of {G}aussian processes and their first-passage
  times.
\newblock {\em Journal of the Royal Statistical Society. Series B
  (Methodological)}, pages 505--522, 1965.

\bibitem{slepian1961first}
D~Slepian.
\newblock First passage time for a particular {G}aussian process.
\newblock {\em The Annals of Mathematical Statistics}, 32(2):610--612, 1961.

\bibitem{ReedSimon}
M~Reed and B~Simon.
\newblock {\em Methods of Modern Mathematical Physics: Scattering theory Vol.
  3}.
\newblock Academic Press, 1979.

\bibitem{NoonZhig}
J~Noonan and A~Zhigljavsky.
\newblock Approximations of the boundary crossing probabilities for the maximum
  of moving sums.
\newblock {\em arXiv preprint arXiv:1810.09229}, 2018.

\bibitem{Quadrature}
JL~Mohamed and LM~Delves.
\newblock {\em Computational Methods for Integral Equations}.
\newblock Cambridge University Press, 1985.

\bibitem{lin1989approximating}
JT~Lin.
\newblock Approximating the normal tail probability and its inverse for use on
  a pocket calculator.
\newblock {\em Applied Statistics}, 38(1):69--70, 1989.

\end{thebibliography}


\begin{thebibliography}{10}

\bibitem{Shepp71}
L.~Shepp.
\newblock First passage time for a particular {G}aussian process.
\newblock {\em The Annals of Mathematical Statistics}, 42(3):946--951, 1971.

\bibitem{slepian1961first}
D.~Slepian.
\newblock First passage time for a particular {G}aussian process.
\newblock {\em The Annals of Mathematical Statistics}, 32(2):610--612, 1961.

\bibitem{li2004lower}
Wenbo~V Li and Qi-Man Shao.
\newblock Lower tail probabilities for {G}aussian processes.
\newblock {\em The Annals of Probability}, 32(1):216--242, 2004.

\bibitem{slepian1962}
D.~Slepian.
\newblock The one-sided barrier problem for {G}aussian noise.
\newblock {\em Bell System Technical Journal}, 41(2):463--501, 1962.

\bibitem{molchan2012survival}
G.~Molchan.
\newblock Survival exponents for some {G}aussian processes.
\newblock {\em International Journal of Stochastic Analysis}, 2012.

\bibitem{landau1970supremum}
H.J. Landau and L.A. Shepp.
\newblock On the supremum of a {G}aussian process.
\newblock {\em Sankhy{\=a}: The Indian Journal of Statistics, Series A},
  32:369--378, 1970.

\bibitem{marcus1972sample}
M.~Marcus and L.A. Shepp.
\newblock Sample behavior of {G}aussian processes.
\newblock {\em In Proc. of the Sixth Berkeley Symposium on Math. Statist. and
  Prob}, 2:423--421, 1972.

\bibitem{adler2007random}
R.J. Adler and J.~Taylor.
\newblock {\em Random Fields and Geometry}.
\newblock Springer, 2007.

\bibitem{pickands1969upcrossing}
J.~Pickands.
\newblock Upcrossing probabilities for stationary {G}aussian processes.
\newblock {\em Transactions of the American Mathematical Society}, 145:51--73,
  1969.

\bibitem{harper2017pickands}
A.J. Harper.
\newblock Pickands’ constant $ {H}_ \alpha$ does not equal $1/{\Gamma
  }(1/\alpha) $, for small $\alpha$.
\newblock {\em Bernoulli}, 23(1):582--602, 2017.

\bibitem{Aldous}
D.~Aldous.
\newblock {\em Probability {A}pproximations via the {P}oisson {C}lumping
  {H}euristic}.
\newblock Springer Science \& Business Media, 1989.

\bibitem{Mehr}
C.B. Mehr and J.A. McFadden.
\newblock Certain properties of {G}aussian processes and their first-passage
  times.
\newblock {\em Journal of the Royal Statistical Society. Series B
  (Methodological)}, 27(3):505--522, 1965.

\bibitem{ReedSimon}
M.~Reed and B.~Simon.
\newblock {\em Methods of Modern Mathematical Physics: Scattering theory Vol.
  3}.
\newblock Academic Press, 1979.

\bibitem{NoonZhig}
J.~Noonan and A.~Zhigljavsky.
\newblock Approximations of the boundary crossing probabilities for the maximum
  of moving sums.
\newblock {\em arXiv preprint arXiv:1810.09229}, 2018.

\bibitem{Quadrature}
J.L. Mohamed and L.M. Delves.
\newblock {\em Computational Methods for Integral Equations}.
\newblock Cambridge University Press, 1985.

\bibitem{beekman1975asymptotic}
J.A. Beekman.
\newblock Asymptotic distributions for the {O}rnstein-{U}hlenbeck process.
\newblock {\em Journal of Applied Probability}, 12(1):107--114, 1975.

\bibitem{pitman2015slepian}
J.~Pitman and W.~Tang.
\newblock The {S}lepian zero set, and {B}rownian bridge embedded in {B}rownian
  motion by a spacetime shift.
\newblock {\em Electronic Journal of Probability}, 20(61):1–28, 2015.

\bibitem{lin1989approximating}
J.T. Lin.
\newblock Approximating the normal tail probability and its inverse for use on
  a pocket calculator.
\newblock {\em Applied Statistics}, 38(1):69--70, 1989.

\end{thebibliography}

\end{document}